\documentclass{birkjour}
\usepackage{amssymb,amsmath}
\usepackage{amsthm}
\usepackage{mathrsfs}
\usepackage{enumerate, indentfirst}
\usepackage[fleqn,tbtags]{mathtools}
\usepackage{color}

\newtheorem*{dfn}{Definition}
\newtheorem*{KvN}{A Krein--von Neumann type extension}

\newtheorem*{MT}{Main Theorem}
\newtheorem*{Cy}{Corollary}

 \theoremstyle{definition}

\theoremstyle{remark}
 
  \numberwithin{equation}{section}

\renewcommand{\phi}{\varphi}
\renewcommand{\theta}{\vartheta}

\DeclareMathOperator{\tform}{\mathfrak{t}}

\DeclareMathOperator{\wform}{\mathfrak{w}}

\DeclareMathOperator{\trace}{Tr}

\DeclarePairedDelimiterX\sipt[2]{(}{)_{\tform}}{#1\,\delimsize\vert\,#2}
\DeclarePairedDelimiterX\sipv[2]{(}{)_{v}}{#1\,\delimsize\vert\,#2}
\DeclarePairedDelimiterX\sipw[2]{(}{)_{w}}{#1\,\delimsize\vert\,#2}

\newcommand{\alg}{\mathscr{A}}

\newcommand{\abs}[1]{\lvert#1\rvert}

\newcommand{\dupC}{\mathbb{C}}

\newcommand{\M}{\mathscr{I}}

\newcommand{\hil}{H}

\newcommand{\bh}{B(\hil)}

\DeclarePairedDelimiterX\sip[2]{(}{)}{#1\,\delimsize\vert\,#2}
\DeclarePairedDelimiterX\siptilde[2]{(}{)_{\!_{\widetilde{A}}}}{#1\,\delimsize\vert\,#2}
\DeclarePairedDelimiterX\sipf[2]{(}{)_{f}}{#1\,\delimsize\vert\,#2}
\DeclarePairedDelimiterX\sipg[2]{(}{)_{g}}{#1\,\delimsize\vert\,#2}
\DeclarePairedDelimiterX\siptw[2]{(}{)_{\tform+\wform}}{#1\,\delimsize\vert\,#2}
\DeclarePairedDelimiterX\set[2]{\{}{\}}{#1\,\delimsize\vert\,#2}
\DeclarePairedDelimiterX\dual[2]{\langle}{\rangle}{#1,#2}
\DeclarePairedDelimiterX\sipa[2]{(}{)_{\!_A}}{#1\,\delimsize\vert\,#2}
\DeclarePairedDelimiterX\sipc[2]{(}{)_{\!_C}}{#1\,\delimsize\vert\,#2}
\DeclarePairedDelimiterX\sipab[2]{(}{)_{\!_{A+B}}}{#1\,\delimsize\vert\,#2}
\DeclarePairedDelimiterX\sipb[2]{(}{)_{\!_B}}{#1\,\delimsize\vert\,#2}

\newcommand{\hsh}{B_2(\hil)}
\newcommand{\trh}{B_1(\hil)}
\newcommand{\finh}{B_F(\hil)}

\allowdisplaybreaks
\begin{document}
\title[A characterization of positive normal functionals on $\bh$]{A characterization of positive normal \\ 
functionals on the full operator algebra}

\author[Z. Sebesty\'en]{Zolt\'an Sebesty\'en}

\address{%
Department of Applied Analysis\\ E\"otv\"os Lor\'and University\\ P\'azm\'any P\'eter s\'et\'any 1/c.\\ Budapest H-1117\\ Hungary}

\email{sebesty@cs.elte.hu}

\author[Zs. Tarcsay]{Zsigmond Tarcsay}
\thanks{ The first author Zsigmond Tarcsay was supported by the Hungarian Ministry of Human
Capacities, NTP-NFTÖ-17. Corresponding author: Tam\'as Titkos.}
\address{%
Department of Applied Analysis\\ E\"otv\"os Lor\'and University\\ P\'azm\'any P\'eter s\'et\'any 1/c.\\ Budapest H-1117\\ Hungary}

\email{tarcsay@cs.elte.hu}

\author[T. Titkos]{Tam\'as Titkos}

\address{%
Alfr\'ed R\'enyi Institute of Mathematics\\
Hungarian Academy of Sciences\\Re\'altanoda utca 13-15.\\
Budapest H-1053\\ Hungary}

\email{titkos.tamas@renyi.mta.hu}



\begin{abstract}
Using the recent theory of Krein--von Neumann extensions for positive functionals we present several simple criteria to decide whether a given positive functional on the full operator algebra $\bh$ is normal. We also characterize those functionals defined on the left ideal of finite rank operators that have a normal extension.
\end{abstract}

\maketitle

The aim of this short note is to present a theoretical application of the generalized Krein--von Neumann extension, namely to offer a characterization of positive normal functionals on the full operator algebra. To begin with, let us fix our notations. Given a complex Hilbert space $\hil$, denote by $\bh$  the full operator algebra, i.e., the $C^*$-algebra of continuous linear operators on $\hil$.
The symbols $\finh,\trh,\hsh$ are referring to the ideals of continuous finite rank operators, trace class operators, and Hilbert--Schmidt operators, respectively. Recall that $\hsh$ is a complete Hilbert algebra with respect to the inner product
\begin{eqnarray*}\label{E:siphsh}
\sip{X}{Y}_2=\trace(Y^*X)=\sum_{e\in \mathcal{E}}\sip{Xe}{Ye},\qquad X,Y\in\hsh.
\end{eqnarray*}
Here $\trace$ refers to the the trace functional and $\mathcal{E}$ is an arbitrary orthonormal basis in $\hil$. Recall also that $\trh$ is a Banach $^*$-algebra  under the norm $\|X\|_1:=\trace (\abs{X})$, and that $\finh$ is dense in both $\trh$ and $\hsh$, with respect to the norms $\|\cdot\|_1$ and $\|\cdot\|_2$, respectively. It is also known that $X\in\trh$ holds if and only if $X$ is the product of two elements of $\hsh$. For the proofs and further basic properties of Hilbert-Schmidt and trace class operators we refer the reader to \cite{KR1}.\\

Let $\alg$ be a von Neumann algebra, that is a strongly closed ${}^*$-subalgebra of $\bh$ containing the identity. A bounded linear functional $f:\alg\to\mathbb{C}$ is called normal if it is continuous in the ultraweak topology, that is $f$ belongs to the predual of $\alg$. It is well known that the predual of $\bh$ is $\trh$, hence every normal functional can be represented by a trace class operator. We will use this property as the definition. 
\begin{dfn}
A linear functional $f:\bh\to\mathbb{C}$ is called a normal functional if there exists a trace class operator $F$ such that
$$f(X):=\trace(XF)=\trace(FX),\qquad X\in\bh.$$
\end{dfn}

\noindent Remark that such a functional is always continuous due to the inequality 
\begin{eqnarray*}\label{trace inequality}
\abs{\trace(XF)}\leq \|F\|_1\cdot\|X\|.
\end{eqnarray*}

Our main tool is a canonical extension theorem for linear functionals which is analogous with the well-known operator extension theorem named after the pioneers of the 20th century operator theory M.G. Krein \cite{Krein} and J. von Neumann \cite{JvN}.  For the details see Section 5 in \cite{SSZT}, especially Theorem 5.6 and the subsequent comments. Let us recall the cited theorem:

\begin{KvN}\label{T:funcext}
    Let $\M$ be a left ideal of the complex Banach $^*$-algebra $\alg$, and consider a linear functional $\varphi:\M\to\dupC$. The following statements are equivalent:
    \begin{enumerate}[\upshape (a)]
    \item There is a representable positive functional $\varphi^{\bullet}:\alg\to\mathbb{C}$ extending $\varphi$, which is minimal in the sense that
          \begin{equation*}
            \varphi^{\bullet}(x^*x)\leq \widetilde{\varphi}(x^*x),\qquad\mbox{holds for all $x\in\alg$,}
          \end{equation*}
          whenever $\widetilde{\varphi}:\alg\to\mathbb{C}$ is a representable extension of $\varphi$.
  
    \item There is a constant $C\geq 0$ such that $|\varphi(a)|^2\leq C\cdot \varphi(a^*a)$ for all $a\in\M$.
    \end{enumerate}

\end{KvN}
We remark that the construction used in the proof of the above theorem is closely related to the one developed in \cite{Sebestyen93} for Hilbert space operators. The main advantage of that construction is that we can compute the values of the smallest extension $\varphi^{\bullet}$ on positive elements, namely
\begin{align}\tag{$*$}\label{F: f_N formula}
\varphi^{\bullet}(x^*x)=\sup\big\{|\varphi(x^*a)|^2\,\big|\,a\in\M,~\varphi(a^*a)\leq 1\big\}\qquad\mbox{for all}~x\in\alg.
\end{align}
The minimal extension $\varphi^{\bullet}$ is called the \emph{Krein--von Neumann extension of $\varphi$}. 

The characterization we are going to prove is stated as follows.
\begin{MT}
For a given positive functional $f:\bh\to\mathbb{C}$ the following statements are equivalent:
\begin{enumerate}[\upshape (i)]
 \item $f$ is normal.
 \item There exists a normal positive functional $g$ such that $f\leq g$.
 \item $f\leq g$ holds for any positive functional $g$ that agrees with $f$ on $\finh$.
 \item For any $X\in\bh$ we have
 \begin{align}\tag{$**$}\label{F: KvN finh}
f(X^*X)=\sup\set{\abs{f(X^*A)}^2}{A\in\finh, f(A^*A)\leq1}. 
\end{align}
\item $f(I)\leq \sup\set{\abs{f(A)}^2}{A\in\finh, f(A^*A)\leq1}$.
\end{enumerate}
\end{MT}
\begin{proof}The proof is divided into three claims, which might be interesting on their own right. Before doing that we make some observations.
For a given trace class operator $S$ let us denote by $f_S$ the normal functional defined by $$f_S(X):=\trace(XS),\qquad X\in\bh.$$
The map\ $S\mapsto f_S$ is order preserving between positive trace class operators and normal positive functionals. Indeed, if $S\geq0$ then
\begin{align*}
f_S(A^*A)=\trace(A^*AS)=\|AS^{1/2}\|^2_2\geq0.
\end{align*}
Conversely, if $f_S$ is a positive functional and $P_{\langle h\rangle}$ denotes the orthogonal projection onto the subspace spanned by $h\in\hil$, we obtain $S\geq0$ by
\begin{align*}
\sip{Sh}{h}=\trace(P_{\langle h\rangle}S)=f_S(P_{\langle h\rangle}^*P_{\langle h\rangle})\geq0,\qquad \mbox{for all}~h\in\hil.
\end{align*}
Our first two claims will prove that (i) and (iv) are equivalent.\\

\textit{\underline{Claim 1.} Let $f$ be a normal positive functional and set $\varphi:=f|_{\finh}$. Then $f$ is the smallest positive extension of $\varphi$, i.e $\varphi^{\bullet}=f$. }

\textit{Proof of Claim 1.} Since $f\geq0$ is normal, there is a positive $S\in\trh$ such that $f=f_S$. By assumption $\varphi$ has a positive extension (namely $f$ itself is one), thus there exists also the Krein--von Neumann extension denoted by $\varphi^{\bullet}$. As $f_S-\phi^{\bullet}$ is a positive functional due to the minimality of $\phi^{\bullet}$, its norm is attained at identity $I$. Therefore it is enough to show that $$\phi^{\bullet}(I)\geq f_S(I)=\trace(S).$$
We know from \eqref{F: f_N formula} that
\begin{equation*}
\phi^{\bullet}(X^*X)=\sup\set{\abs{\phi(X^*A)}^2}{A\in\finh, \phi(A^*A)\leq1}
\end{equation*}
for any $X\in\bh$. Choosing $A=\trace(S)^{-1/2}P$ for any   projection $P$ with finite rank, we see that $\phi(A^*A)=\trace(S)^{-1}\trace(PS)\leq1$, whence
\begin{align*}
\phi^{\bullet}(I)\geq \abs{\phi(A)}^2=\frac{\trace(PS)^2}{\trace(S)}.
\end{align*}
Taking supremum in $P$ on the right hand side we obtain $\phi^{\bullet}(I)\geq\trace(S)$, which proves the claim.\\

\textit{\underline{Claim 2.} The smallest positive extension of $\varphi$, i.e. $(f|_{\finh})^{\bullet}$ is normal.}

\textit{Proof of Claim 2.} First observe that the restriction of $f$ to $\hsh$ defines a continuous linear functional on $\hsh$ with respect to the norm $\|\cdot\|_2$. Due to the Riesz representation theorem, there exists a unique representing operator $S\in\hsh$ such that
\begin{align}\tag{$***$}\label{E:fS}
f(A)=\sip{A}{S}_2=\trace(S^*A),\qquad \textrm{for all $A\in\hsh$}.
\end{align}
We are going to show that $S\in\trh$. Indeed, let $\mathcal{E}$ be an orthonormal basis in $\hil$ and let $\mathcal{F}$ be any non-empty finite subset of $\mathcal{E}$. Denoting by $P_{\mathcal{F}}$ the orthogonal projection onto the subspace spanned by $\mathcal{F}$ we get
\begin{align*}
 \sum_{e\in\mathcal{F}} \sip{Se}{e}=\sip{P_{\mathcal{F}}}{S}_2=f(P_{\mathcal{F}})\leq f(I).
\end{align*}
Taking supremum in $\mathcal{F}$ we obtain that $S$ is in trace class. By Claim 1, the smallest positive extension $\phi^{\bullet}$ of $\phi$ equals $f_S$ which is normal. This proves Claim 2.

Now, we are going to prove (ii)$\Rightarrow$(i).\\

\textit{\underline{Claim 3.} If there exists a normal positive functional $g$ such that $f\leq g$ holds, then $f$ is normal as well.}

\textit{Proof of Claim 3.}  Let $g$ be a normal positive functional dominating $f$, and let $T$ be a trace class operator such that $g=f_T$. According to Claim 2 it is enough to prove that $f=\phi^{\bullet}$. Since $h:=f-\phi^{\bullet}$ is positive, this will follow by showing that $h(I)=0$. We see from \eqref{E:fS} that  $h(A)=0$ for any finite rank operator $A$. Consequently, as $h\leq f\leq f_T$, it follows that 
\begin{align*}
h(I)=h(I-P)\leq f_T(I-P)=\trace(T)-\trace(TP),
\end{align*}
for any finite rank  projection $P$. Taking infimum in $P$ we obtain $h(I)=0$ and therefore Claim 3 is established.\\

Completing the proof we mention all the missing trivial implications. Taking $g:=f$, (i) implies (ii). As \eqref{F: KvN finh} means that $\varphi^{\bullet}=f$, equivalence of (iii) and (iv) follows from the minimality of the Krein-von Neumann extension. Replacing $X$ with $I$ in \eqref{F: KvN finh} we obtain (v). Conversely, (v) implies (iv) as $\varphi^{\bullet}\leq f$ and $f-\varphi^{\bullet}$ attains its norm at $I$. \end{proof}

Finally, we remark that the above proof contains a characterization of having normal extension for a functional defined on $\finh$.
\begin{Cy}
Let $\varphi:\finh\to\mathbb{C}$ be a linear functional. The following statements are equivalent to the existence of a normal extension.
\begin{enumerate}[\upshape (a)]
\item There is a $C\geq0$ such that $|\varphi(A)|^2\leq C\cdot \varphi(A^*A)$ for all $A\in\finh$.
\item There is a positive  functional $f$ such that $f|_{\finh}=\varphi$.
\item There is an $F\in\trh$ such that $\varphi(A)=\trace(FA)$ for all $A\in\finh$.
\end{enumerate}
\end{Cy}

\bibliographystyle{abbrv}

\begin{thebibliography}{10}




\bibitem{KR1}
R. V. Kadison and J. R. Ringrose,
\newblock  {\em Fundamentals of the theory of operator algebras I.},
\newblock  Academic Press, New York, 1983.

\bibitem{Krein} M. G. Krein, \emph{The theory of self-adjoint extensions of semi-bounded Hermitian transformations
and its applications, I-II}, Mat. Sbornik 20, 431--495, Mat. Sbornik 21, 365--404 (1947) (Russian)

\bibitem{JvN}
J. von Neumann, 
\newblock Allgemeine Eigenwerttheorie Hermitescher Funktionaloperatoren, \newblock \emph{Math. Ann.}, \textbf{102} (1930) 49--131.





\bibitem{Sebestyen93}
Z.~Sebesty{\'e}n,
\newblock Operator extensions on {H}ilbert space,
\newblock {\em Acta Sci. Math. (Szeged)}, \textbf{57} (1993), 233--248.

\bibitem{SSZT}
Z.~Sebesty{\'e}n, Zs. Sz\H ucs, and Zs. Tarcsay,
\newblock Extensions of positive operators and functionals,
\newblock \emph{Linear Algebra Appl.}, \textbf{472} (2015), 54--80.

\end{thebibliography}

\end{document}